\documentclass[11pt]{article}
\usepackage[utf8]{inputenc}
\usepackage{amsthm}
\usepackage{amsmath}
\usepackage{amssymb}
\usepackage{enumitem}
\usepackage{bookmark}
\usepackage[small]{titlesec}
\RequirePackage{geometry}
\usepackage{caption}
\usepackage{graphicx}
\usepackage{hyperref}
\usepackage[T1]{fontenc}
\usepackage{amsfonts}
\usepackage{braket}
\usepackage[l2tabu, orthodox]{nag}
\usepackage{pdf14}
\usepackage{xpatch}
\RequirePackage{tabto}
\usepackage{indentfirst}
\usepackage{amscd,stmaryrd,latexsym}
\usepackage{xcolor}

\newtheorem{thm}{Theorem}
\newtheorem{lem}{Lemma}[section]

\newtheorem{cor}{Corollary}[section]
\newtheorem{Prop}{Proposition}[section]
\newtheorem{Prop*}{Proposition}
\newtheorem{Def}{Definition}[section]
\newtheorem*{theorem*}{Theorem}
\newtheorem*{corollary*}{Corollary}
\theoremstyle{definition}
\newtheorem{rem}{Remark}[section]

\newtheorem{eg}{Example}[section]

\newcommand{\inlinesubsection}[1]{%
  \refstepcounter{subsection}% αυξάνει τον μετρητή
  \addcontentsline{toc}{subsection}{\protect\numberline{\thesubsection}#1}% ToC
  \textbf{\thesubsection\ #1. }% εμφανίζεται στην ίδια γραμμή
}

\newcommand\res{\mathop{\hbox{\vrule height 7pt width .5pt depth 0pt
\vrule height .5pt width 6pt depth 0pt}}\nolimits}

\DeclareMathOperator{\Per}{\mathrm{Per}}
\DeclareMathOperator{\spt}{\mathrm{spt}}
\DeclareMathOperator{\sing}{\mathrm{sing}}
\DeclareMathOperator{\reg}{\mathrm{reg}}

\DeclareMathAlphabet{\mathscr}{OT1}{pzc}{m}{it}

\begin{document}

\title{\textbf{Smooth approximation of one--sided minimizing CMC hypersurfaces with isolated singularities}}
\author{ Konstantinos Leskas\footnote{Department of Mathematics, National and Kapodistrian University of Athens, Greece.}}
\date{}

\maketitle

\begin{abstract} 
\noindent 
We consider a CMC hypersurface with an isolated singularity, given by the boundary of a set of finite perimeter, that minimizes the prescribed mean curvature functional on one side. If a varifold tangent at the singularity is a non--flat regular cone of multiplicity--1 that minimizes perimeter on the side occupied by the blow--up, then we show that, in a small neighborhood of the singular point, there exists a sequence of smooth CMC hypersurfaces, with the same prescribed mean curvature, that converges to the initial hypersurface. We also establish the full CMC variational equation and regularity for minimizers of a class of singular obstacle problems.
    
\end{abstract}

\section{Introduction}

The formation of singularities in variational constructions is a central difficulty in the field of geometric analysis. In particular, it is known that hypersurfaces of $\mathbb R ^{n+1}$ that minimize the perimeter functional, or more generally the prescribed mean curvature functional, have a singular set of Hausdorff dimension at most $n-7$. From the work of Bombieri, De Giorgi and Giusti \cite{BDG}, the bound on the size of the singular set is sharp. A natural question is whether such singularities can be removed by small perturbations of the boundary data or the ambient metric. 

Hardt and Simon proved in \cite{HS} that, for a $7$--dimensional hypersurface of $\mathbb R ^8$ with an isolated singularity that minimizes area and has a $6$--dimensional smooth boundary $\Gamma $, a small perturbation of $\Gamma$ yields an area minimizer that is smooth. This result was later used by Smale \cite{Smale} to establish  generic regularity for area minimizers in non--zero integral homology classes of closed, oriented, $8$--dimensional manifolds, by perturbing the ambient metric. The generic regularity theorem for area minimizing hypersurfaces has since been established in ambient dimensions $9$ and $10$ by Chodosh, Mantoulidis and Schulze \cite{CMS} and in ambient dimension $11$ by Chodosh, Mantoulidis, Schulze and Wang \cite{CMSW}. An extension of the Hardt--Simon perturbation result has been proven by Liu \cite{Liu} for stationary one--sided area minimizers. 

For CMC hypersurfaces, Bellettini and the author \cite{BL} proved smooth approximation near isolated singularities of minimizers of the prescribed mean curvature functional. Such approximations are useful in Bellettini and Marshall--Stevens' \cite{BK} generic existence theorem for smooth CMC hypersurfaces in $8$--dimensional closed manifolds with positive Ricci curvature. Here, we address the smooth approximation problem under the weaker assumption that the hypersurface minimizes only on one side.  

\medskip 
We set up some notation, see also Section \ref{preli} for more details. For a bounded open set $U$ and $\lambda \in \mathbb R$, we consider the prescribed mean curvature functional in the class of sets of finite perimeter $A \subset U,$
\[J_\lambda (A) := \Per_U(A) - \lambda |A|.\]  
We say that the set $A$ minimizes $J_\lambda$ on one side if it attains the infimum of $J_\lambda$ among all competitor sets of finite perimeter $F\subset B_2$, such that $F\setminus (A\cap B_1) = A\setminus (A\cap B_1).$ Let $\nu_A$ be the inward--pointing unit normal vector field of a set of finite perimeter $A$ and $\partial ^* A$ its reduced boundary. Let $V_{\mathcal A}$ be the varifold associated to the free part of the boundary of $A$ that we get after removing $\partial^* A \cap(B_2 \setminus \overline B_1).$ We say that $V_{\mathcal A} $ is a critical point of $J_\lambda$ if it satisfies the CMC variational equation \eqref{variational eq}.

\begin{thm}\label{main theorem}Let $A \subset \mathbb R^8$ be a set of locally finite perimeter and $\lambda \in \mathbb R.$ Assume that, near the origin, $\partial A$ is smoothly embedded, with an isolated singularity at the origin, $V_{\mathcal A}$ is a critical point of $J_\lambda$ and $A$ minimizes $J_\lambda$ on one side. Suppose that a varifold tangent at the origin is given by a non--flat regular cone $C$ of multiplicity--1 that locally minimizes the perimeter functional on the side occupied by the blow--up of $A.$ 

Then there exists $r>0$ and sets of locally finite perimeter $A_j \subset A$ that minimize $J_\lambda$ on one side in $B_r$, with $\chi_{A_j\cap B_r} \xrightarrow[L^1]{}\chi_{A\cap B_r}$ and $\partial A_j \cap B_r$ are smoothly embedded CMC hypersurfaces with the same prescribed mean curvature $\lambda.$ Moreover, $\partial A_j $ converges smoothly to $\partial A$ on any compact subset of $B_r \setminus \{0\}.$ 
\end{thm} 

From the construction of CMC hypersurfaces with isolated singularities in \cite{Irv}, which extends the work of \cite{CHS}, together with the minimizing argument in \cite{Lesk} and \cite{HS}, there are examples of CMC hypersurfaces that satisfy the assumptions of Theorem \ref{main theorem}. These examples are asymptotic to the Lawson cone $C_{1,5}=\{(x,y)\in \mathbb R^2\times \mathbb R^6 \mid 5|x|^2=|y|^2\}$, which is known to minimize area on one side, see \cite{Dav}. Since $C_{1,5}$ is not two--sided area minimizing, these examples lie outside the minimizing setting of \cite{BL}.

We do not establish a general multiplicity--1 theorem for critical one--sided minimizers of $J_\lambda.$ In our proof, the tangent cone assumption is used to establish a weaker multiplicity--1 convergence result, Theorem \ref{mult--1}, which is then used in the proof of the singular maximum principle, Corollary \ref{singular max pri}.

In particular, if $A$ is a one--sided minimizer of $J_\lambda$ and  $V_{\mathcal A}$ is a critical point of $J_\lambda$, there could be rescalings $A_j = r_j^{-1}A$, $r_j \searrow 0,$ for which 
\[\chi_{A_j} \xrightarrow[L^1_{\text{loc}}]{} 1 , \hspace{2mm} \partial\llbracket A_j\rrbracket \rightharpoonup 0, \hspace{2mm} V_{\mathcal A_{j}}\to 2|C|,\]
where $|C|$ denotes a non--flat minimal regular cone. Thus, oppositely oriented sheets would cancel in the limit as currents while contributing a multiplicity--2 varifold mass. At a touching singularity, see Section \ref{2.2}, any rescaling of $A$ satisfies this behavior but the limit $|C|$ is flat. 

For $\lambda=0$, such a behavior is not possible for stationary one--sided area minimizers, due to \cite[Theorem B]{Ilm}. However, for $\lambda \neq 0$, a possible mechanism would be two CMC sheets on opposite sides of an area minimizing cone $C$ meeting only at the vertex, and $A$ is the complement of the region between the two sheets. Their scalar mean curvatures with respect to the induced orientation from the cone $C$ would have opposite signs, hence it does not contradict the singular maximum principle, Corollary \ref{singular max pri}. Constructing such a configuration, with the required one--sided minimality for $A$, remains a separate issue.

\medskip 
\noindent
\textit{Question:} Must every varifold tangent cone of a critical one--sided minimizer of $J_\lambda$ at an isolated non--touching singularity have multiplicity--1? 

\medskip
In order to prove Theorem \ref{main theorem}, we need to investigate the behavior of the first variation of minimizers with respect to a singular obstacle. Unlike the situation of smooth obstacle problems, see for example \cite{Duzaar}, \cite{BL}, the singularities of the obstacle raise extra difficulties. The main result that we establish, which does not require any multiplicity--1 assumption on the varifold tangent cone, is the following:

\begin{thm}\label{singular obstacle problems} Let $A\subset B_2 \subset \mathbb R^8$ be a one--sided minimizer of $J_\lambda$, with $0<|\lambda|<7.$ Assume that $\partial A$ is smoothly embedded with an isolated, non--touching, singularity at the origin and meets $\partial B_1$ transversely. Assume also that $V_{\mathcal A}$ is a critical point of $J_\lambda.$ For any set of finite perimeter $E_0 \subset A$ there exists a minimizer $E\subset A $ that minimizes $J_\lambda$ among sets $F\subset B_2$ with $F\setminus (A\cap B_1) = E_0 \setminus (A\cap B_1)$. 

Let $V_{\mathcal E}$ be the associated varifold to the free part of the boundary of $E$ and $T_b$ its fixed current boundary, as defined in Section \ref{preli}. Then 
\[\delta V_{\mathcal E}(X) = - \lambda\int X\cdot \nu_E \: d\mu_{V_{\mathcal{E}}}, \]
for any $X \in C^1_c(B_2 \setminus \spt T_b ; \mathbb R ^8).$ Moreover, $\spt V_{\mathcal E}\setminus \spt T_b$ is smoothly embedded except at isolated singular points. 
    
\end{thm}

Theorem \ref{singular obstacle problems} is proven in Section \ref{2.3}. We show that the variational equation holds across the isolated singularity of $\partial A$ at the origin using a method similar to Simon's \cite{Simon unfin.}. A barrier argument excludes contact with the corner of the obstacle formed by $\partial A \cap \partial B_1$. 

Once Theorem \ref{singular obstacle problems} is established we then show a multiplicity--1 convergence result, Theorem \ref{mult--1}, using a comparison argument similar to Liu's \cite{Liu}. We can then follow the same methods of \cite{HS}, \cite{BL}, to show an instance of a singular maximum principle, Corollary \ref{singular max pri}. In proving the smooth approximation theorem, Theorem \ref{main theorem}, we perturb the initial smooth boundary data given by $\partial A \cap \partial B_1$ on one side and construct minimizers with this boundary data and obstacle given by $\overline {A\cap B_1}.$ Theorem \ref{singular obstacle problems} and the singular maximum principle exclude contact at the origin. Then we repeat the blow--up argument of \cite{HS}, \cite{BL}. The extended Hardt--Simon foliation result, proven by Liu \cite{Liu}, together with Allard's interior regularity theorem \cite{Allard1}, exclude singular points from accumulating at the origin. 

\medskip 
Although the proof is presented in $\mathbb R^8$, since it is essentially a result about isolated singularities, Theorem \ref{main theorem} extends to $\mathbb R ^{n+1}$ for $n\geq 7$, under the same hypotheses of an isolated singularity and a varifold tangent at the origin that is a non--flat regular cone of multiplicity--1 that minimizes perimeter on one side; see Remark \ref{any dimension}.

\medskip
\noindent
\textbf{AI Statement.} ChatGPT assisted with the exposition and revision of this manuscript and suggested the auxiliary vector field used in the proof of Lemma \ref{avoidance}, which was then checked by the author. The remaining arguments and proofs were developed by the author who takes full responsibility for the final manuscript and its mathematical content. 

\medskip
\noindent 
\textbf{Acknowledgments.} The material in this work partly overlaps with the thesis of the author who would like to thank Costante Bellettini for first introducing him to this problem during his PhD years at University College London. The author would also like to thank Panagiotis Gianniotis for several stimulating discussions during the preparation of this work.

\section{Singular obstacle problems}

\noindent 
In this section we define the prescribed mean curvature functional $J_\lambda$ and introduce the notion of a one--sided minimizer. We then investigate the existence and regularity of minimizers for certain obstacle problems. 

\medskip

\noindent
\inlinesubsection{Preliminaries}\label{preli} Let $U$ be an open subset of $\mathbb R ^{n+1}$ and $E \subset U$ a Lebesgue measurable set. We say that $E$ is a set of locally finite perimeter in $U$ if the characteristic function $\chi_E \in BV_{\text{loc}}(U),$ where $BV_{\text{loc}}(U)$ denotes the space of functions of locally bounded variation. If $\chi_E \in BV(U)$ we say that $E$ is a set of finite perimeter.

By Riesz's representation theorem, if $E$ is a set of locally finite perimeter in $U$, then there exist a Radon measure $\mu_E$ and a $\mu_E$--measurable vector field $\nu_E$, with $|\nu_E|=1$ $\mu_E$--a.e. in $U$, such that 

\begin{equation}
\label{int by parts}
\int_E \text{div} X \: d\mathcal L^{n+1} = -\int_ U \nu_E \cdot X \: d\mu_E,
\end{equation}
for any vector field $X \in C_c^1(U ; \mathbb R ^{n+1})$, where $\mathcal L ^{n+1}$ denotes the Lebesgue measure in $\mathbb R ^{n+1}.$ The vector field $\nu_E$ is called the inward--pointing  unit normal vector field of $E.$ We define the relative perimeter of $E$ in $K$ as $\Per_K(E):= \mu_E(K),$ where $K\subset U.$ 

Let $E_j$ be a sequence of sets of locally finite perimeter. We say that $E_j \xrightarrow[j \to \infty]{} E$ if and only if $\chi_{E_j}\xrightarrow[L^1_{\text{loc}}]{}\chi_E.$ In particular, we have the standard compactness for sets of locally finite perimeter under uniform perimeter bounds (see \cite[Theorem 12.26, Corollary 12.27]{Maggi}).

We define the prescribed mean curvature functional $J_\lambda$, for $\lambda \geq 0$, as 

\[ J_\lambda (E) := \Per_{U}(E) - \lambda |E|,\]
where $E \subset U$ is a set of finite perimeter and $|E|:= \mathcal L^{n+1}(E). $

Let $B_r := B^{n+1}_r(0)$ denote the $(n+1)$--dimensional ball of $\mathbb R^{n+1}$ of radius $r$ centered at the origin. We are interested in the following minimization problem 

\begin{equation}
\label{min problem}
\inf \bigg \{ J_\lambda (E) \mid E \subset B_2, \: \chi_{E} |_{B_2 \setminus (A\cap B_1)} = \chi_{E_0} |_ { B_2\setminus (A \cap B_1)} \bigg \},
\end{equation}
where $A \subset B_2$ and  $E_0 \subset A$ are sets of finite perimeter. We note that the equality between the characteristic functions is understood in the $\mathcal L^{n+1}$--a.e. sense. In general, when we write that two sets of finite perimeter are equal, we will mean that the characteristic functions are the same in the $\mathcal{L}^{n+1}$--a.e. sense. 

Standard compactness and the lower semi--continuity of the perimeter functional yield the following existence result (see \cite[Lemma 2.1]{BL} for more details).

\begin{Prop}\label{existence} There exists a set of finite perimeter $E \subset B_2$ such that $E\setminus (A\cap B_1) = E_0 \setminus (A \cap B_1)$ and $J_\lambda(E)$ equals the infimum in \eqref{min problem}.  

\end{Prop}

To each set of locally finite perimeter $E$ in $U$, we associate its reduced boundary $\partial^* E$ in $U$; see \cite[Section 15]{Maggi}. Then from \cite[Proposition 12.19]{Maggi}, we may always assume that $\spt \mu_E\cap U = \overline{\partial ^* E}\cap U=\partial E\cap U$. More importantly, from De Giorgi's structure theorem \cite[Theorem 15.9]{Maggi}, we have that $\partial^* E$ is a countably $n$--rectifiable set with multiplicity--1  and $\mu_E = \mathcal H^n \res\partial^* E,$ where $\mathcal H^n$ denotes the Hausdorff measure of dimension $n$ in $\mathbb R^{n+1}$. In particular, the perimeter functional is $\Per_K(E) = \mathcal H^n(\partial^* E \cap K),$ for any $K \subset U$.

Thus, to each set of locally finite perimeter $E$, we can associate an integer rectifiable current $T_E := (\partial ^* E, 1, -\star \nu_E ),$ where we use the notation of \cite[Chapter 6.3]{Simon's notes} ($\star$ denotes the Hodge star so that $\nu_E \wedge \star \nu_E$ gives the positive orientation of $\mathbb R^{n+1}$).  We can also associate an integer rectifiable varifold $V_E := \underbar v(\partial ^* E, 1), $ where we use the notation of \cite[Chapter 4]{Simon's notes}. 

Furthermore, if $F$ is a competitor of finite perimeter for the minimization problem \eqref{min problem}, we set  $\mathcal F := \partial ^* F \setminus [\partial^* E_0 \cap (B_2 \setminus \overline{A\cap B_1})]$. We consider the associated integer rectifiable varifold $V_{\mathcal F}:= \underbar v(\mathcal F, 1)$ and the associated integer rectifiable current $T_{\mathcal F} := (\mathcal{F}, 1, -\star \nu_F)$. Note that  $\partial^* F \cap( B_2\setminus \overline{A\cap B_1}) = \partial^* E_0 \cap (B_2 \setminus \overline{A\cap B_1}),$ where equality is understood in the $\mathcal H^n$--a.e. sense. In particular,  
\[T_{\mathcal F} = \partial \llbracket F\rrbracket - \partial \llbracket E_0\rrbracket \res ( B_2\setminus \overline{A\cap B_1}),\]
where $\llbracket F \rrbracket$ denotes the current given by Lebesgue integration over $F$. Note that $T_{\mathcal F}$  has a current boundary given by $T_b:= - \partial (\partial \llbracket E_0\rrbracket\res(B_2\setminus \overline {A\cap B_1})).$ Thus, for a competitor set $F$ and a minimizer $E$ of \eqref{min problem}, using the above notation, the minimization condition becomes 
\[\mathbb M(T_{\mathcal E}) - \lambda |E\cap B_1|\leq \mathbb M(T_{\mathcal F}) - \lambda |F \cap B_1|,\]
where $\mathbb M(\cdot)$ denotes the mass of a current. 

\medskip
\noindent
\inlinesubsection{One--sided minimizers and their properties}\label{2.2}We say that a set of finite perimeter $A\subset B_2$ is a one--sided minimizer of $J_\lambda$ if 
\[J_\lambda(A) \leq J_\lambda (E),\]
for any set of finite perimeter $E \subset B_2$ with $A\setminus (A\cap B_1) = E\setminus (A\cap B_1).$  

We turn our attention to the variational equations satisfied by a one--sided minimizer. Let $V_{\mathcal A}, T_{\mathcal A}$ and $T_b$ be as defined in Section \ref{preli}. 

\begin{Def} We say that $V_{\mathcal A}$ is a critical point of $J_\lambda$ if for any $X \in C^1_c(B_2 \setminus \spt T_b; \mathbb R ^{n+1})$ 
\[\delta V_{\mathcal A}(X)= - \lambda \int \nu_A \cdot X \:d\mu_{V_{\mathcal A}},\]
where $\delta V_{\mathcal A}(\cdot)$ denotes the first variation of the varifold $V_{\mathcal A}.$ If $\lambda=0$, we say that $V_{\mathcal A}$ is stationary. 
\end{Def}

In particular, from the well--known first variation formula \cite{Simon's notes}, we have that for a critical point $V_{\mathcal A}$ 
\begin{equation}
\label{variational eq}
 \int \text{div}_{\partial^* A} X \: d \mu_{V_{\mathcal A}} +\lambda \int \nu_{A}\cdot X \:d\mu_{V_{\mathcal A}}=0,
\end{equation}
for any $X \in C^1_c(B_2 \setminus \spt T_b; \mathbb R ^{n+1}).$

It is not necessarily true that for a one--sided minimizer $A$ the varifold $V_{\mathcal A}$ satisfies the variational equations \eqref{variational eq}.

\begin{eg}\label{example}Consider the union of two half--lines $l_1: y=ax$ and $l_2: y=-ax$, $a>0$. Let $L$ be the non--convex set with boundary the union of the half--lines $l_1, l_2$. Then 
\[\Per_{B_2}(L) \leq \Per_{B_2}(E),\]
for any set $E \subset B_2$ with $E\setminus (L\cap B_1)=(L\cap B_2)\setminus (L\cap B_1).$  In particular, $L\cap B_2$ is a one--sided minimizer for $\lambda=0.$ However, the varifold $V_{\mathcal L}$ does not satisfy \eqref{variational eq}. Indeed, if we take a suitable vector field $X$ that is compactly supported around the origin, then 
\[\delta V_{\mathcal L}(X) = -X(0)\cdot (n_1+n_2)\neq 0, \]
where $n_1, n_2$ are the inward unit co--normals at the origin. 

\medskip
We can similarly construct an example with $\lambda>0. $ Take $\lambda = \frac{1}{2}$ and consider a set $A\subset B_2\subset \mathbb R ^2$ so that $A\cap \overline B_1$ corresponds to the set given by $\overline B_1 \cap [(\overline B_2(c^+) \cap \{y\geq 0\})\cup (\overline B_2(c^-) \cap \{y\leq 0\})]$, where $c^\pm = (1, \pm \sqrt 3)$ and $A$ has no additional boundary pieces in $\partial B_1$ coming from $A \cap (B_2 \setminus \overline B_1).$ In particular, $\partial A$ in $\overline B_1$ is the union of two arcs in $\overline B_1$ that meet at the origin. We show that $A$ minimizes $J_\lambda$ on one side.  

Since the arcs are graphical, we choose the following vector field defined in the ball $B_{1+\epsilon}$
\[X(x,y) = 
\begin{cases}
    \bigg( - \dfrac{\sqrt{4-(y+\sqrt 3)^2}}{2}, \dfrac{y+\sqrt{3}}{2} \bigg) & \text{ if } y<0 \\
    \bigg( -\dfrac{\sqrt{4-(y-\sqrt{3})^2}}{2}, \dfrac{y-\sqrt{3}}{2} \bigg) &\text{ if }  y \geq 0.
\end{cases}\]
Note that the vector field $X$ is singular along the symmetry line given by the $x$--axis. We have that $|X|_\infty =1$ and $X= - \nu_A$ along the two arcs. We compute its divergence in a distributional sense. In particular, we have that as a Radon measure 
\[d(\text{div}X) = \lambda \: d \mathcal  L^2\res B_{1+\epsilon}-\sqrt{3}\:d\mathcal{H}^1 \res (\{y=0\}\cap B_{1+\epsilon})\leq \lambda \:d \mathcal  L^2. \]
Thus, for the vector field $X$ we may apply the Gauss--Green formula of \cite{Crasta}. In particular, for a competitor set $F\subset B_2$ such that $F \setminus (A\cap B_1) = A\setminus (A\cap B_1)$, we get that 

\[\int_{\overline B_1} \chi_A\: d(\text{div}X) +\int_{\overline B_1}(X, D\chi_A) = \int_{\overline B_1} \chi_F \:d(\text{div}X) + \int_{\overline B_1} (X, D\chi_F),\]
where $(X, D\chi_A)$ and $(X, D\chi_F)$ denote the pairing defined in \cite{Anzel}; see also \cite{Crasta}. 

Using the definition of the pairing and the vector field $X$, we get that $(X, D\chi_A) = - |D\chi_A|$ and $(X, D\chi_F) \geq - |D\chi_F|.$ Putting these together, we get that
\[\int_{\overline B_1} (\chi_A - \chi_F)\:d(\text{div}X)\geq \Per_{\overline B_1}(A)- \Per_{\overline B_1}(F).\]
Since $d(\text{div}X) \leq \lambda \: d\mathcal L ^2$ and $F\subset A$, we get that 
\[\Per_{\overline B_1}(F) - \lambda |F \cap B_1| \geq \Per_{\overline B_1}(A)-\lambda |A\cap B_1|.\]
Adding and subtracting the terms in $B_2 \setminus \overline B_1$ where the two sets agree, we get that 
\[J_\lambda(F)\geq J_\lambda (A),\]
which shows the one--sided minimizing property of $A.$ 

Note that at the origin the associated varifold $V_\mathcal A$ does not satisfy the variational equations \eqref{variational eq}. If it did, then blowing up at the origin would give a stationary limit. However, the limit is the union of two half--lines, which we have already shown is not stationary at the origin.

\end{eg}

We now recall certain types of singularities. We say that $p \in \spt V$ is a classical singularity of a varifold $V$ if $\spt V \cap B^{n+1}_r(p)$, for some $r>0$, is the union of three or more embedded $C^{1, \alpha}$--hypersurfaces with common boundary $L$ such that $p \in L$ and at least two of them meet transversely along $L$.

We say that $p \in \spt V$ is a touching singularity of the varifold $V$ if $\spt V \cap B^{n+1}_r(p)$ is the union of two embedded $C^{1, \alpha}$-- hypersurfaces $M_1, M_2$ that meet tangentially at $p$ and $(\overline {M_i} \setminus M_i) \cap B^{n+1}_r(p)=\emptyset,$ $i=1,2.$ We denote by $\sing_T V$ the set of touching singularities. %and by $\genreg V := \reg V\cup \sing_T V$ the generalized regular set.

We can have one--sided minimizers of $J_\lambda$ for $\lambda>0$ with touching singularities. 

\begin{eg}Let $A \subset B_2 \subset \mathbb R ^2$ be such that $A \cap \overline B_1 = \overline B_1 \cap [\overline B_2 (c^+)\cup \overline B_2(c^-)]$, where $c^\pm = (\pm 2, 0)$ and $A$ has no additional boundary pieces in $\overline B_1$ coming from $A\cap(B_2 \setminus \overline B_1)$. In particular, $\partial A$ in $\overline B_1$ is the union of the two arcs given by the disks $\overline B_2 (c^\pm)$ that meet at the origin tangentially. We show that $A$ minimizes $J_\lambda$ on one side for $\lambda = \frac{1}{2}$. 

We define the following vector field 
\[X(x,y) = \begin{cases}
    \bigg( - \dfrac{\sqrt{4-y^2}}{2}, \dfrac{y}{2} \bigg) & \text{ if } x\geq 0 \\
    \bigg( \dfrac{\sqrt{4-y^2}}{2}, \dfrac{y}{2} \bigg) &\text{ if }  x < 0.
\end{cases}\]

As in Example \ref{example}, we treat $X$ as a singular vector field and we show it satisfies all the properties of a singular sub--calibration. Then the exact same computation, using the Gauss--Green theorem of \cite{Crasta}, shows the minimizing property. Indeed, we have that as a Radon measure the divergence satisfies
\[d(\text{div}X) = \lambda\: d\mathcal L^2 - \sqrt{4-y^2}\:d\mathcal H^1\res \{x=0\}\leq \lambda\: d\mathcal L^2.\]
Furthermore, $X= - \nu_A$ along the two arcs. This is enough to repeat the same argument of Example \ref{example}.

\end{eg}

We now show that a one--sided minimizer does not have classical singularities. 

\begin{lem}\label{absence of classical} Let $V_{\mathcal A}$ be the associated varifold to a one--sided minimizer $A$. Then $\sing_C V_{\mathcal A}=\emptyset.$ 
\end{lem}

\begin{proof} Let $A_t$, $t\geq 0$, be a one--sided deformation of $A$ that agrees with $A$ outside $B_1$ and such that the map $ t \to \chi_{A_t}$ is continuous in the $L^1_{\text{loc}}$--sense. Then $J_\lambda(A_t) \geq J_\lambda (A)$, and if the right derivative of the perimeter functional exists, we get that 
\begin{equation}
\label{L^1-stationary}
\dfrac{d}{dt}\bigg|_{t=0^+} J_\lambda (A_t) \geq 0.
\end{equation}

In particular, we may construct a deformation that is admissible and breaks apart the singularity as in \cite[Section 9]{BW1}. This will contradict the variational inequality. We briefly explain the computation. 

Assume $p\in \sing_C V_{\mathcal A}$ and choose coordinates so that the common boundary has tangent plane $x_1=x_2=0$ and two of the hypersurfaces $M_1, M_2$ have tangent planes given by the half--hyperplanes $H_{\pm}=\{x_2=\pm ax_1 \mid x_1\geq 0\}$ where $a>0$ and $L=\{x_2=x_1=0\}$. Let $B$ be a small ball around $p$ and $E$ a connected component of $A \cap B$ given by the hypersurfaces $M_1, M_2$. Assume that no other hypersurface with boundary exists in $B \cap \{-ax_1<x_2<ax_1 \}.$ 

Consider the vector field $X = \Psi \partial_{x_1}$, where $\Psi \in C^1_c(B)$, $\Psi\geq 0$ and $\Psi\equiv 1$ in a neighborhood of $p.$ Define $A_t = (A \setminus E)\cup ((Id + t\Psi \partial_{x_1})(E))$, for $t\geq 0.$ Then $A_t$ is a one--sided deformation satisfying the boundary condition, since outside $E$ it agrees with $A$. Thus, \eqref{L^1-stationary} is satisfied for $A_t$.  In particular, for $E_t = (Id + t\Psi \partial_{x_1})(E)$ we still have \eqref{L^1-stationary}.

Blowing up $E$ around $p$, we get a limit cone $C$ given by $H_\pm$. Due to the weak convergence of the first variations as functionals of $C^1_c(B)$ combined with \eqref{L^1-stationary}, we get that $\delta C(X)\geq 0$. The latter is not possible since, as in Example \ref{example}, the first variation is 
\[\delta C (X) = - \dfrac{2}{\sqrt{1+a^2}}\int_L X \cdot e_1 \:d\mathcal H ^{n-1}\] 
which, for $X=\Psi \partial _{x_1}$ is negative.
\end{proof}

Let now $V_{\mathcal A}$ be the associated varifold of a one--sided minimizer that satisfies \eqref{variational eq}. We assume that there are no touching singularities for $V_{\mathcal A}$. We say that $V_{\mathcal A}$ is stable if for $\Sigma:=\reg V_{\mathcal A}\setminus \spt T_b$ we have that

\begin{equation} 
\label{stab ineq}
\int_{\Sigma} |A_\Sigma|^2\eta^2 \:d\mathcal H ^n \leq \int_\Sigma | \nabla^\Sigma \eta|^2 \:d \mathcal H^n,
\end{equation}
where $A_\Sigma$ denotes the second fundamental form of $\Sigma$, $\nabla^\Sigma $ denotes the Levi--Civita connection of $\Sigma$ and $\eta \in C^{1}_c(\Sigma).$ 

\begin{rem}\label{rem on stability} Let $\nu $ be the inward--pointing unit normal vector field that orients $\Sigma$. Then \eqref{stab ineq} is equivalent to the non--negativity of the second variation of $J_\lambda$ with respect to variations $X= \eta \nu$, where $\eta \in C^1_c(\Sigma)$. See \cite{Bar} for more details. Note that a standard approximation argument shows that \eqref{stab ineq} holds for $\eta \in C^{0,1}_c(\Sigma).$
\end{rem}

\begin{lem}
\label{stability} Let $V_{\mathcal A}$ be the varifold associated to a one--sided minimizer $A$. Assume that $V_{\mathcal A}$ has no touching singularities and satisfies \eqref{variational eq}. Then $ V_{\mathcal A}$ is stable. 
\end{lem}

\begin{proof} We first show that \eqref{stab ineq} is valid for variations $X = \eta \nu$, where $\nu$ is the inward--pointing unit normal vector field and $\eta \in C^1_c(\Sigma), \eta \geq 0.$ Indeed, if $\theta_t$ is the flow of $X$ then $\spt (\theta_{t *} V_{\mathcal A}) \subset \overline{A \cap B_1}$. Thus, if we glue the boundary data, we obtain sets of finite perimeter $A_t$ such that $A_t \setminus (A\cap B_1) = A \setminus (A\cap B_1)$. Hence 
\[J_\lambda (A_t) \geq J_\lambda (A)\]
and 
\[\frac{d}{dt}\bigg|_{t=0}J_\lambda (A_t)=0.  \] 
The above equations imply that 
\[\frac{d^2}{dt^2}\bigg|_{t=0} J_\lambda (A_t) \geq 0.\]
From Remark \ref{rem on stability}, the latter is equivalent to \eqref{stab ineq} for any $\eta \in C^{0,1}_c(\Sigma)$ with $\eta\geq 0$. 

For an arbitrary $\eta \in C^{0,1}_c(\Sigma)$, we write $\eta=\eta^+-\eta^-$, where $\eta^+,\eta^- \geq 0$, $\eta^+, \eta^- \in C^{0,1}_c(\Sigma)$, $\eta^+\eta^- \equiv 0 $ and $\nabla^\Sigma \eta^+ \cdot \nabla^\Sigma \eta^- \equiv 0$. Then 
\[\int_\Sigma |\nabla^\Sigma\eta|^2-\int_\Sigma |A_\Sigma|^2 \eta^2 = \]
\[=\int_\Sigma |\nabla^\Sigma \eta^+|^2 +\int_\Sigma |\nabla^\Sigma \eta^-|^2 -\int_\Sigma |A_\Sigma|^2 (\eta^+)^2-\int_\Sigma |A_\Sigma|^2(\eta^-)^2 \geq 0, \]
where we used the fact that \eqref{stab ineq} is valid for $\eta^+, \eta^-$ as shown above. \end{proof}

From Lemma \ref{absence of classical} and Lemma \ref{stability}, we may apply the regularity theory \cite[Theorem 2.2]{BW1} to get the following:

\begin{thm}\label{reg theorem 1} If $V_{\mathcal A}$ is the associated varifold of a one--sided minimizer $A$ that has no touching singularities and satisfies \eqref{variational eq}, then $\spt V_{\mathcal A}\setminus  \spt T_b$ is smooth except on a set $S$ of Hausdorff dimension at most $n-7.$ If $n=7$, the set $S$ consists of isolated singular points.
\end{thm}

\noindent 
\inlinesubsection{Singular obstacle problems}\label{2.3} In what follows, let $A$ be a one--sided minimizer of $J_\lambda$ for $0<\lambda < n$ whose associated varifold $V_{\mathcal A}$ satisfies \eqref{variational eq}. Assume that $\partial A $ is smoothly embedded except at the origin, which is not a touching singularity, and that $\partial A \cap \partial B_1$ is a transverse intersection.

From Proposition \ref{existence}, there exists a set of finite perimeter $E$ attaining the infimum in \eqref{min problem}.  Let $V_{\mathcal E} $ be the associated varifold  and $T_{\mathcal E}$ the associated current  with current boundary $T_b$ as defined in Section \ref{preli}. 

The minimization property implies that for any $X \in C_c^1(B_2 \setminus \spt T_b; \mathbb R ^{n+1})$  such that $\theta_{t *}V_{\mathcal E} \subset \overline {A\cap B_1}$, where $\theta_t$ is the flow of $X$ and $t\geq 0,$ 

\begin{equation}
\label{sign of 1st var}
    \delta V_{\mathcal E}(X) \geq - \lambda \int \nu_E \cdot X \:d\mu_{V_{\mathcal E}}.
\end{equation}
We want to extend the variational inequality \eqref{sign of 1st var} so that we have an equality. 

First, note that if $X \in C_c^1 (A\cap B_1 ; \mathbb R ^{n+1})$,  then $\phi_t(p)=p+t X(p)$ is a two--sided deformation for $t$ small enough. Hence we have equality in \eqref{sign of 1st var}. We next analyze the cases in which the variation field moves a point $p \in \partial B_1 \cap A$ or  $p \in \partial A  \cap B_1$ with $p \neq 0$.

\begin{lem}\label{var part I} Let $X \in C^1_c(B_2 \setminus \spt T_b \setminus (\partial A \cap \partial B_1) \setminus \{0\}; \mathbb R ^{n+1})$. Then we have that 
\[\delta V_{\mathcal E} (X)= -\lambda \int \nu_E \cdot X \:d\mu_{V_{\mathcal E}}.\] 
    \end{lem}

\begin{proof} Choose a point $p \in \partial B_1 \cap A \cap  \spt V_{\mathcal E} \setminus \spt T_b$ and a small ball $B_r(p) \subset \subset A$ which does not intersect $\spt T_b.$ Then, for any $X \in C_c^1 (B_r(p) ; \mathbb R ^{n+1})$  such that $X \cdot \nu_{B_1}\geq 0$, we have that $\spt (\theta_{t *} V_{\mathcal E}) \subset \overline {A\cap B_1}$, so \eqref{sign of 1st var} holds.

We can then apply \cite[Theorem 7]{White} to $\partial B_1$ and $\spt V_{\mathcal E}$ in $B_r(p)$. This implies that $\partial B_1$ must have constant mean curvature $\lambda$ in $B_r(p)$, contradicting the assumption that $\lambda < n.$  

Now choose a point $p \in \partial A \cap B_1 \cap \spt V_{\mathcal E}\setminus \spt T_b$ and a small ball $B_r(p) \subset \subset B_1$ that does not intersect $\spt T_b$ and does not contain the origin. Then for any $X \in C_c^1 (B_r(p); \mathbb R ^{n+1})$ such that $X \cdot \nu_A \geq 0$ we have that $\spt (\theta_{t *} V_{\mathcal E}) \subset \overline {A\cap B_1} $ so \eqref{sign of 1st var} holds. 

In $B_r(p)$, we are in a smooth obstacle setting similar to that of \cite[Lemma 2.2]{BL} or \cite[Lemma 3.3]{Duzaar}. Applying the same arguments, based on the variational inequality and Riesz's representation theorem, we obtain
\[\delta V_{\mathcal E}(X) + \lambda \int X \cdot \nu_E \:d\mu_{V_{\mathcal E}} = \int X \cdot \nu_A \: d\mu,\]
where $\mu$ is a Radon measure with $\mu \leq (\text{div}_{\mathcal E}\nu_A +\lambda \nu_E \cdot \nu_A)\:d \mathcal H ^n\res \mathcal E \cap \partial A.$ 

At $\mathcal H^n$--a.e. point in $\mathcal E \cap \partial A$ the tangent spaces agree and, since $E \subset A$, we have $\nu_E=\nu_A$. Consequently, $\text{div}_{\mathcal E}\nu_A = \text{div}_{\partial A}\nu_A = -\lambda.$ Hence $\mu\equiv 0$, which implies \eqref{variational eq} for such vector fields. 
\end{proof}

Next, we extend the monotonicity formula for $ V_{\mathcal E}$ around the origin. Our proof is analogous to the one in \cite{Simon unfin.}, where the case $\lambda=0$ is treated. 

\begin{lem}\label{var eq part II} For any $\tau<\rho<1$, we have that 
\[e^{\lambda\tau}\tau^{-n}\mu_{V_{\mathcal E}}(B_\tau)\leq e^{\lambda \rho}\rho^{-n}\mu_{V_{\mathcal E}}(B_\rho).\]
    
\end{lem}

\begin{proof} Let $\rho_0 <1$ so that $B_{\rho_0}$ does not intersect $\spt T_b.$ For $0< \sigma < \rho<\rho_0$, choose a function $\phi_\sigma$ such that $\phi_\sigma \equiv 0$ for $t< \frac{\sigma}{2}$, $\phi_\sigma \equiv 1$ for $\sigma < t$,  $\phi_\sigma ' \geq 0$ and a function $\psi$ such that $\psi\equiv 0$ for $t\geq 1$, $\psi \equiv 1$ for $t \in [0, 1-\epsilon)$ and $\psi'\leq 0.$  Let $\psi_{\sigma, \rho}(t)=\phi_\sigma (t) \psi(\frac{t}{\rho})$ and define the vector field $X(x) = \psi_{\sigma, \rho}(|x|)x.$ Then the vector field is supported in $B_\rho \setminus \overline B_{\frac{\sigma}{2}}$ and so $X \in C^1_c(B_2 \setminus \spt T_b \setminus (\partial A \cap \partial B_1) \setminus \{0\}; \mathbb R ^{n+1})$. Thus, by Lemma \ref{var part I}, we have that \eqref{variational eq} holds for $X$. In particular, 

\[\int \text{div}_{\mathcal E} X \: d\mu_{V_{\mathcal{E}}} +\lambda \int X \cdot \nu_{E} \:d\mu_{V_{\mathcal{E}}}=0.\]
Setting $\overrightarrow H= \lambda \nu_{E}$ and computing the divergence, we get that 
\[\int n \psi_{\sigma, \rho} (|x|)\: d\mu_{V_{\mathcal{E}}} +\int \psi'_{\sigma, \rho}(|x|)\cdot |x| \cdot |\nabla^{\mathcal E}|x||^2\:d\mu_{V_{\mathcal{E}}}+\]
\[+\int \psi_{\sigma, \rho}(|x|) x \cdot \overrightarrow H \:d\mu_{V_{\mathcal{E}}}=0. \]
After multiplying by $\rho^{-n-1}$ and simplifying we obtain 
\[-\frac{d}{d\rho}\bigg(\rho^{-n}\int \psi_{\sigma, \rho}(|x|)\:d\mu_{V_{\mathcal{E}}}\bigg)+\rho^{-n}\frac{d}{d\rho}\int\psi_{\sigma, \rho}(|x|)(1-|\nabla^{\mathcal E}|x||^2)\:d\mu_{V_{\mathcal{E}}}+\] 
\[+ E_0(\rho)\leq 0,\] 
where $\displaystyle E_0(\rho)=\rho^{-n}\int \psi_{\sigma, \rho}(|x|)\:\dfrac{x}{\rho}\cdot \overrightarrow H \:d\mu_{V_{\mathcal{E}}}.$ Since 
\[\rho^{-n}\frac{d}{d\rho}\int\psi_{\sigma, \rho}(|x|)(1-|\nabla^{\mathcal E}|x||^2)\:d\mu_{V_{\mathcal{E}}}\geq 0,\] 
we get that 
\[\frac{d}{d\rho}\bigg(\rho^{-n}\int \psi_{\sigma, \rho}(|x|)\:d\mu_{V_{\mathcal{E}}}\bigg)\geq  E_0(\rho) \geq -\lambda \rho^{-n} \int \psi_{\sigma, \rho}(|x|) \:d\mu_{V_{\mathcal E}}.\]
Thus 
\[\frac{d}{d\rho}\bigg(e^{\lambda \rho}\rho^{-n}\int \psi_{\sigma, \rho}(|x|)\:d\mu_{V_{\mathcal{E}}}\bigg)\geq 0.\]

For $\sigma< \tau <\rho$, we integrate the latter inequality from $\tau $ to $\rho$ and use the monotonicity of $\psi$ to obtain
\[e^{\lambda\rho}\rho^{-n}\int \psi_{\sigma, \rho}(|x|)\:d\mu_{V_{\mathcal{E}}} \geq e^{\lambda \tau} \tau^{-n} \int \psi_{\sigma, \tau}(|x|)\: d\mu_{V_{\mathcal{E}}}.  \]
As $\sigma \to 0$, we have that $\phi_\sigma \to 1$ for $\mathcal H^n$--a.e. point in $\mathcal E$. Since $\mathcal E$ has finite $\mathcal H^n$--measure, we use the dominated convergence theorem to get that 
\[e^{\lambda \rho}\rho^{-n} \int_{\mathcal E}\phi_\sigma (|x|)\psi \bigg(\frac{|x|}{\rho}\bigg)\:d\mathcal H^n \xrightarrow[\sigma \to 0]{}e^{\lambda \rho}\rho^{-n}\int_\mathcal E \psi \bigg(\frac{|x|}{\rho}\bigg)\:d\mathcal H^n.   \] 
If we let $\epsilon \to 0$, then $\psi \bigg(\dfrac{|x|}{\rho}\bigg) \to \chi_{B_\rho}$. Since $\mathcal E$ has finite $\mathcal H^n$--measure, we use the dominated convergence theorem again to get that 
\[e^{\lambda \rho}\rho^{-n} \mathcal H^n (\mathcal E \cap B_\rho)\geq \tau^{-n}e^{\lambda \tau} \mathcal H^n(\mathcal E \cap B_\tau),\]
which is what we wanted. 
\end{proof}

\noindent
Putting together Lemma \ref{var part I} and Lemma \ref{var eq part II}, we get the following: 

\begin{lem}\label{final var eq} For any $X \in C^1_c(B_2 \setminus \spt T_b \setminus (\partial A \cap \partial B_1) ; \mathbb R ^{n+1})$, we have that \[\delta V_{\mathcal E}(X) = -\lambda \int \nu_E \cdot X\: d\mu_{V_{\mathcal E}}.\] 
    \end{lem}

\begin{proof} We only need to consider the case in which the vector field $X$ moves the origin, so we may assume that $\spt X \subset B_\rho$ for some sufficiently small $\rho<1$.  We choose a function $\eta_\epsilon \in C^1(B_\rho)$ such that $0\leq \eta_\epsilon \leq 1$, $\eta_\epsilon \equiv 0$ in $B_{\epsilon/2}$, $\eta_\epsilon \equiv 1$ in $B_\rho\setminus B_\epsilon$ and $|\nabla \eta_\epsilon |\leq \frac{C}{\epsilon}.$ 
Thus the vector field $\tilde X = \eta_\epsilon X $ is admissible and from Lemma \ref{var part I} we have that 
\[\int_\mathcal E \eta_\epsilon\text{div}_\mathcal E X \:d\mathcal H^n + \int_\mathcal E \nabla^\mathcal E \eta_\epsilon \cdot X\:d\mathcal H^n +\lambda \int_\mathcal E \eta_\epsilon X \cdot \nu_{E}\:d\mathcal H^n =0.\] 
For the second term we use Lemma \ref{var eq part II}, to get that 
\begin{equation}
\label{cutting the singular set}
\int_\mathcal{E} |\nabla^\mathcal E \eta_\epsilon |\cdot |X| \:d\mathcal H^n \leq \frac{C}{\epsilon}|X|_{\infty} \mathcal H^n(\mathcal E \cap B _\epsilon) \leq C \epsilon^{n-1}.  
\end{equation}
As $\epsilon \to 0$, we have $\eta_\epsilon \to 1$ for $\mathcal H^n$-a.e. point. Hence the dominated convergence theorem and \eqref{cutting the singular set} give the result. 
\end{proof} 

Finally, in order to show that \eqref{variational eq} holds for vector fields that move points in $\partial A \cap \partial B_1$, we will show that $\spt T_{\mathcal E}$, away from its boundary, avoids the singular part of the obstacle given by $\partial A \cap \partial B_1$. In order to do that, we construct an appropriate barrier to $\partial A \cap \partial B_1$ whose normal points inside the corner. The proof is inspired by the proof of \cite[Theorem 2]{White}.  

\begin{lem}(Avoidance lemma)\label{avoidance}  $(\spt V_{\mathcal E}\setminus \spt T_b)\cap \partial A \cap \partial B_1 = \emptyset$.
    
\end{lem}

\begin{proof} Let $p \in (\partial A \cap \partial B_1)\setminus \spt T_b$. In a neighborhood $U$ of $p$ disjoint from $\spt T_b$, the signed distance functions $d_1$ of $\partial A$ and $d_{2}$ of $\partial B_1$ are smooth. In particular, the obstacle $\overline {A \cap B_1}$ in $U$ is given by the set $\Omega = \{x\in U \mid d_1(x)\geq 0 \text{ and } d_2(x)\geq 0\}.$ 

After a translation, we may assume that $p=0.$ By transversality, $\nabla d_1 (0)$and $ \nabla d_2 (0)$ are linearly independent. We define the following orthonormal vectors 
\[v= \dfrac{\nabla d_1(0) +\nabla d_2(0)}{|\nabla d_1(0) +\nabla d_2(0)|}\]
\[w=\dfrac{\nabla d_2(0) - \nabla d_1(0)}{|\nabla d_2(0)-\nabla d_1(0)|}.\]
Then we may write 
\[\nabla d_2(0) = a v + bw\]
\[\nabla d_1(0) = av - bw\]
where $a , b>0.$ 

For any $x \in \Omega$, Taylor expansion of $d_1, d_2$ at $0$ gives
\[\nabla d_1(0) \cdot x \geq - O(|x|^2)\]
\[\nabla d_2(0) \cdot x \geq - O(|x|^2).\]
From these inequalities and the expressions of $\nabla d_1 (0)$ and $ \nabla d_2(0)$ in terms of $v$ and $ w$, we conclude that 
\begin{equation}
\label{dont know}
v\cdot x \geq \frac{b}{a}|w\cdot x|- C_0 |x|^2.
\end{equation}
We now define the following polynomial 
\[u(x) = x\cdot v + L |x|^2 - M (x\cdot w)^2,\]
where $L = C_0+1$ and $M$ will be chosen later. 
From \eqref{dont know}, we have that 
\[u(x) \geq \frac{b}{a}|w\cdot x| - M(x\cdot w)^2 + |x|^2.\] 
If $x \in \Omega \cap B_r$ and $r \leq \frac{b}{2aM}$, then we get that 
\begin{equation}
    \label{localization} 
    u(x) \geq |x|^2
    \text{, hence } u(x)=0 \text{ if and only if } x=0. 
\end{equation}

Since $\nabla u = v + 2L x - 2M(x\cdot w)w$, we get that $\nabla u(0) = v$ and so $\nabla d_1(0) \cdot \nabla u(0)=a >0.$  Similarly, $\nabla d_2(0) \cdot v >0.$ Thus from continuity, choosing a smaller $r$ if necessary, we get that 
\begin{align}
    \label{push to 1 side}
    &\nabla d_1(x) \cdot \nabla u(x) >0\\
    &\nabla d_2(x) \cdot \nabla u(x) >0
\end{align}
for all $x \in B_r.$

The Hessian of $u$ is given by 
\[\nabla ^2 u = 2L\: Id - 2M w\otimes w.\]
We want to compute the trace $tr_ S (\nabla ^2 u)$ with respect to any $n$--dimensional plane $S.$ Let $(e_i)$ be an orthonormal basis of $S$. Then 
\[tr_S (\nabla ^2 u) = \sum_{i=1}^n <\nabla^2 u(e_i), e_i>= 2Ln - 2M \sum_{i=1}^n (w\cdot e_i)^2=\]
\[= 2Ln - 2M |P_S w|^2.\]
If $\xi$ is a unit normal to $S$, then we rewrite 
\[tr_S (\nabla ^2 u) = 2Ln - 2M + 2M (\xi \cdot w)^2.\]
Since $v$ and $w$ are orthonormal and $\xi$ is a unit vector, we have $(w\cdot \xi)^2 +(v\cdot \xi)^2\leq 1$ and hence $|P_S v|^2 \geq (w\cdot \xi)^2.$ 

We define the function $F : B_r \times Gr(n, n+1) \to \mathbb R$ by 
\[F(x, S) = tr_S(\nabla ^2 u(x)) - 2M |P_S \nabla u(x)|^2 + \lambda |\nabla u(x)|.\] 
Note that the function is uniformly continuous. At $x=0$, for any plane $S$, we have 
\[F(0,S) = tr_S (\nabla ^2 u(0)) - 2M |P_S v|^2 + \lambda = \]
\[=2nL -2M + 2M (\xi \cdot w)^2 - 2M |P_S v|^2 +\lambda \leq   \]
\[\leq 2nL +\lambda - 2M < -2,\]
where we chose $2M > 2nL +\lambda +2>0.$ 
Thus from uniform continuity and a sufficiently small $r$, we get that $F(x, S) \leq -1$ for any $x \in B_r$ and any $n$--plane $S$.

We define the vector field $X(x) = \phi(x) f(u(x)) \nabla u(x)$, where $\phi \in C^\infty _c (B_r)$, $0 \leq \phi \leq 1$ and $\phi \equiv 1$ in $B_{\frac{r}{2}}$ and 
\[f(t) = \begin{cases}
e^{-2M t - \frac{1}{\epsilon-t}+\frac{1}{\epsilon}} & \text{ if } t<\epsilon \\
0 & \text{ if } t\geq \epsilon
\end{cases}\]
where $\epsilon < \frac{r^2}{4}.$ Thus $f(0)=1$ and $f'(t)\leq - 2M f(t).$ 

From \eqref{push to 1 side}, we have that the flow of $X$ pushes points into the obstacle $\overline {A\cap B_1}$ which means that \eqref{sign of 1st var} is valid for $X.$ Note also that in $\Omega \cap (B_r\setminus B_{\frac{r}{2}})$ we have from \eqref{localization} that $u\geq \frac{r^2}{4}>\epsilon$, hence $f(u) = f'(u)=0.$ Thus we compute for any $n$--plane $S$ and any $x\in B_{\frac{r}{2}}$ 
\[\text{div}_S X + \lambda |X| = f(u)tr_S (\nabla ^2 u) + f'(u) |P_S \nabla u(x)|^2 + \lambda |\nabla u(x)|f(u)\leq  \]
\[\leq f(u) tr_S (\nabla ^2 u) - 2M f(u) |P_S \nabla u|^2 + \lambda f(u) |\nabla u(x)|= f(u) F(x, S) \leq - f(u).\]  
Thus we get that 
\begin{equation}
    \label{curvature estimate}
    \text{div}_S X + \lambda |X| \leq - f(u).
\end{equation}
If we now use $X$ in \eqref{sign of 1st var} together with \eqref{curvature estimate} for the varifold $V_{\mathcal E}$, we get that 
\[\int_{B_r}f(u) \:d\mu_{V_{\mathcal E}}\leq 0.\] 
Since $f(0)=1$, this implies that $\mu_{V_{\mathcal E}}(B_\delta)=0$, for some small $\delta< r$, which means that $0 \notin \spt V_{\mathcal E}.$
\end{proof}

If we combine Lemma \ref{var part I}, Lemma \ref{final var eq} and Lemma \ref{avoidance}, we obtain that \eqref{variational eq} holds for any $X \in C^1_c(B_2 \setminus \spt T_b; \mathbb R ^{n+1})$, hence proving the first part of Theorem \ref{singular obstacle problems} in any dimension.

\begin{thm} \label{reg theorem 2} We have that $\sing_T V_{\mathcal E}\setminus \spt T_b = \emptyset$ and $\spt V_{\mathcal E} \setminus \spt T_b$ is smooth except on a set $S$ of Hausdorff dimension at most $n-7.$ If $n=7$, the set $S$ consists of isolated singular points. 
\end{thm}

\begin{proof} Once we show the absence of touching singularities for $V_{\mathcal E}$, the result follows from Theorem \ref{reg theorem 1}. Let $p \in \sing_T V_{\mathcal E} \cap A\cap B_1$. Locally, we have a two--sided minimizer hence we cannot have touching singularities due to the density estimates \cite[Theorem 21.11]{Maggi}. 

If $p \in \sing_T V_{\mathcal E} \cap \partial A\cap B_1 \setminus \{0\}$, then the minimization property for $E$ implies that the $(n+1)$--density of $E$ is $1$ at $p.$ Since $p \in \partial A$, the density of $A$ at $p$ is equal to $\frac{1}{2}$, and since $E \subset A$ we get a contradiction.  

If $p=0\in \sing_T V_{\mathcal E}$, then the same density argument gives that $E$ has density $1$ at the origin. Blowing up at the origin and using $E \subset A$, we get that the limit varifold tangent of $\partial A$ at the origin is a plane, possibly with multiplicity. From the sheeting theorem and higher regularity theorem \cite[Theorem 3.1, Theorem 3.3]{BW1} applied to $V_{\mathcal A}$, we have that the origin is a regular point or a touching singularity, contrary to our assumption. 
\end{proof}

\section{Multiplicity--1 tangent cones} 
\label{3}
\noindent
Let $A$ be a one--sided minimizer of $J_\lambda$. Let $V_{\mathcal A}$ and $T_{\mathcal A}$ be as defined in Section \ref{preli}. Assume that $V_{\mathcal A}$ satisfies \eqref{variational eq}, $\sing_T V_{\mathcal A}=\emptyset$ and $0 \in \sing V_{\mathcal A}.$ We restrict our attention to the case that $n=7$ and so the ambient space is $\mathbb R ^8$. The reason is that from Theorem \ref{reg theorem 1} the singular set consists of isolated points. Furthermore, we assume that a varifold tangent $C$ of $V_{\mathcal A}$ at the origin has multiplicity--1. From the compactness theorem of \cite{BW1} and regularity theory of \cite{BW1} and \cite{Wick2}, it follows that $C$ is a regular cone (smooth everywhere except at the origin). From standard algebraic topology, see \cite{BL}, we have that $\mathbb R^8\setminus C$ has two connected components $D_\pm$. Thus for the associated current (still denoted by $C$), we have that $C= \partial \llbracket D \rrbracket$ where $D= D_+ \subset \mathbb R ^8$ and the current tangent is the same as the varifold tangent. From \cite{Simon unique tangents} we conclude that the varifold tangent $C$ is unique. 

We prove that a minimizer $E$ of the minimization problem \eqref{min problem} with obstacle $A$ satisfying the above properties has a multiplicity--1 varifold tangent at the origin given by $C$ (provided that $0 \in \overline{\partial^*E}$). 

\begin{thm}\label{mult--1}(Multiplicity--1 tangent cones) Let $E_j := \eta_{0, \rho_j}(E)$ be a sequence of rescalings, where $\rho_j \searrow 0$. Then, after passing to a subsequence,  we have that 
\[T_{\mathcal{E}_j} \rightharpoonup \partial \llbracket  D\rrbracket,  \]
\[V_{\mathcal{E}_j} \to C.\] 
\end{thm}

\begin{proof}  From Theorem \ref{reg theorem 2} and since the result is local, we assume that the origin is the only isolated singularity of $V_{\mathcal E}$ (else we restrict to a small neighborhood of the origin).

From Lemma \ref{final var eq} we have that \eqref{variational eq} holds around the origin and so from the monotonicity of the mass ratio we have that 
\[\sup_j\Per_{B_R}(E_j) <+\infty, \hspace{2mm} \text{ for any } R>0.\]
Thus we can apply the standard compactness for sets of finite perimeter \cite[Corollary 12.27]{Maggi} to get that 
\begin{align*}
&\chi_{E_j} \xrightarrow[L^1_{loc}]{}\chi_{E_\infty}\\
&T_{\mathcal{E}_j} \rightharpoonup T_\infty = \partial \llbracket E_\infty \rrbracket. 
\end{align*}

By Theorem \ref{reg theorem 2}, $V_{\mathcal E_j}$ has no touching singularities. Together with Lemma \ref{absence of classical} and Lemma \ref{stability}, this allows us to apply the compactness theorem \cite[Theorem 2.3]{BW1}. Thus $V_{\mathcal E _j} \to V_{\infty}$, where the varifold $V_\infty$ is stationary, stable and has no classical singularities. Then the regularity theorem of \cite{Wick2} implies that the singular set of $V_{\infty}$ consists of isolated singular points. Due to the monotonicity formula of the mass ratio around the origin, Lemma \ref{var eq part II}, we know that $V_{\infty}$ is a cone. Since the singular set of $V_{\infty}$ consists of isolated points and is a cone, we conclude that $\spt V_\infty$ is a regular cone. 

Since $\spt V_{\mathcal E _j} \to \spt V_\infty$ in the Hausdorff distance sense and $E_j \subset \eta_{0,\rho_j}(A)$, the uniqueness of the varifold tangent of $V_{\mathcal A}$ implies that $\spt V_\infty \subset \overline D.$ In particular, we can apply the maximum principle for stationary varifolds with a small singular set, \cite{Ilm}, to get that $\spt V_\infty = C.$ By the constancy theorem of \cite{Simon's notes}, we conclude that  $V_{\infty} = k C,$ where $k \in \mathbb Z _{\geq 1}.$  

We show that $E_\infty = D.$ From the minimization property we have the standard uniform lower bounds on volume, see \cite[Theorem 21.11]{Maggi}. Thus $E_\infty \neq \emptyset.$ From $E_\infty \subset D$ we have that $E_\infty \neq \mathbb R ^8.$ From \cite{White2} we have that $\spt T_\infty \subset \spt C$. Thus the distributional derivative $D\chi_{E_\infty}=0$ in $\mathbb R ^8 \setminus C = D \cup D^c.$ By the constancy theorem for functions of locally bounded variation, we have that $\chi_{E_\infty}$ is constant on each connected component of $\mathbb R ^8 \setminus C.$ Since $E_\infty \neq \emptyset$ and  $E_\infty \neq \mathbb R^8$, we conclude that $\chi_{E_\infty}=1$ in $D$ and $\chi_{E_\infty}=0$ in $D^c$, and hence $E_\infty = D$.

We now show that $k=1.$  From the Hausdorff distance convergence we have that for any compact $K' \subset \mathbb R ^{8}$
\begin{equation}\label{hausd}
    \spt V_{\mathcal E_j}\cap K' \subset \{|d| < \epsilon_j\}\cap K'  
\end{equation}
where $d$ denotes the signed distance of $C$, positive on $D$, and $\epsilon_j \to 0.$

Consider $K'= \overline B_{2r}\setminus B_{ \delta/2} $ and $K= \overline B_r \setminus B_\delta$ so that the intersection of $\spt V_{\mathcal E_j}$ with $\partial B_r, \partial B_{ \delta}$ is a smoothly embedded $(n-1)$--submanifold (from Sard's theorem and the smoothness of $\spt V_{\mathcal E_j}$ away from the origin we can find such a uniform $r$). The signed distance function $d$ of $C$ is smooth in a tubular neighborhood $U$ of $C$ in the annulus $B_{2r} \setminus B_{ \delta/2}$. 

Define $F_j:= \{d>2\epsilon_j \}\cap K'$. We show that $F_j \subset E_j\cap K'$. Note that the sets $\{d>\epsilon_j\}\cap K'$ and $\{d< - 2\epsilon_j\}\cap K'$ are connected since the cone is regular. From \eqref{hausd}, it follows that $D\chi_{E_j}=0$ in $\{d> \epsilon_j\}\cap K'.$ By the constancy theorem for BV--functions, we conclude that $\chi_{E_j}$ is constant and since $\chi_{E_j} \to \chi_D$ we conclude that $\chi_{E_j}=1$ in $\{d>\epsilon_j\}\cap K'$. Similarly, we have that $\chi_{E_j} =0$ in $\{d<-2\epsilon_j \}\cap K'.$ This shows the desired inclusion and that $(E_j \setminus F_j)\cap K' \subset \{|d|< 2 \epsilon_j\}\cap K'.$

We estimate the perimeter of $F_j$ in the interior $\operatorname{int}(K)$ of $K$. Notice that $\partial^* F_j \cap \operatorname{int}(K) = \Phi_j(C\cap K_j)$, where $\Phi_j(x) = x+ 2\epsilon_j \nu_C (x)$ and $K_j $ is the annulus $B_{\sqrt{r^2-4\epsilon_j^2}}\setminus B_{\sqrt{\delta^2-4\epsilon_j^2}}.$ Using the area formula for the perimeter functional, we have that 
\[\Per_{\operatorname{int}(K)}(F_j) = \int_{C \cap K_j}J_C \Phi_j \:d \mathcal{H}^7\]
and the Jacobian is given by 
\[\text{det}(Id+2\epsilon_j \nabla_C \nu_C)=1+2\epsilon_j H_C+O(\epsilon_j^2).\] 
Since $C$ is stationary, we conclude that 
\[\Per_{\operatorname{int}(K)} (F_j)\leq \Per_{K_j} (D) + O(\epsilon_j^2).\]

Define the set $R_j = (E_j\setminus K)\cup (F_j \cap K)$. Then $R_j$ is a competitor for $J_{\lambda_j}$, where $\lambda_j = \rho_j \lambda$. Thus 
\[J_{\lambda_j}(E_j) \leq J_{\lambda_j}(R_j) \text{ so }\]
\[\Per_K(E_j) -\lambda_j |E_j\cap K| \leq \Per_K (R_j)-\lambda_j |R_j \cap K|.\]
We may further assume that $\partial^* F_j$ intersects $\partial K$ transversely. Using the perimeter formulas of \cite[Theorem 16.3]{Maggi}, the fact that $F_j \subset E_j$, and the transversality assumptions, we get that 
\[\Per_K (E_j) - \lambda_j |E_j \cap K|\leq \Per_{\operatorname{int}(K)}(F_j)+ \mathcal H ^7 (\partial K \cap (E_j \setminus F_j)) -\lambda_j |F_j \cap K|.\]

Since $(E_j \setminus F_j )\cap \partial K \subset \{|d|<2\epsilon_j\}\cap \partial K$ and the volume terms tend to zero, letting $j\to \infty$ gives
\[k \Per_K(D)=\lim_j \Per_K(E_j) \leq \lim_j \Per_{\operatorname{int}(K)} (F_j) \leq\]\[\leq \lim_j (\Per_{K_j} (D) + O(\epsilon_j ^2)) =  \Per_K (D).  \] 
This shows that $k=1$ and finishes the proof. \end{proof}

\begin{cor}(singular maximum principle)\label{singular max pri} For $0<\lambda<7$, let $E$ be a minimizer of $J_\lambda$ with respect to an obstacle $A$ satisfying the properties above. Assume that the boundaries of $T_{\mathcal A}$ and $T_{\mathcal E}$ are smoothly embedded, connected $(n-1)$--submanifolds. If  $0 \in \overline{\partial^* E}\cap \overline{\partial^* A}$, then  $A\cap B_1 = E\cap B_1.$ 
\end{cor}

\begin{proof} \emph{Step 1.} We show that $\spt T_{\mathcal A} $ and $\spt T_{\mathcal E}$ are connected. Assume to the contrary that $\spt T_{\mathcal A} = K_1\cup K_2$, where $K_1, K_2$ are non--empty, disjoint, compact sets with $d(K_1, K_2)>0.$ Since the boundary is connected, we may assume that $\spt \partial T_{\mathcal A} \subset K_1.$ Define the vector field $X(x)= \phi(x) x,$ where $\phi \in C^1_c(B_2)$, with $\phi\equiv 0$ in a neighborhood of $K_1$ and $\phi\equiv 1$ in a neighborhood of $K_2.$ In particular, the vector field is compactly supported away from the current boundary of $T_{\mathcal A}$.  We can apply the variational equations \eqref{variational eq} for $X$ to get that  
\[7\mu_{V_\mathcal A}(K_2)=\int \text{div}_{\partial A} X\:d\mu_{V_{\mathcal A}}=-\lambda \int X \cdot \nu_A \:d\mu_{V_\mathcal A} \leq \lambda \mu_{V_\mathcal A}(K_2).\]
Since $\lambda <7$, we get that $\mu_{V_\mathcal A}(K_2)=0$, which is not possible unless $\spt V_{\mathcal A}$ is connected.

From Lemma \ref{avoidance}, we have that $\eqref{variational eq}$ is valid for $V_{\mathcal E}$. We can use the exact same argument to show that $\spt V_{\mathcal E}$ is connected. 

\medskip
\noindent
\emph{Step 2.} From Theorem \ref{mult--1}, we repeat the same arguments of \cite[Proposition 4.1]{BL} to get that there exists $\rho>0$ such that $\spt T_{\mathcal A}\cap B_\rho = \spt T_{\mathcal E}\cap B_\rho.$ For $\epsilon>0$, define the set 
\[\mathcal Z_{\epsilon} = \bigg\{x \in \spt T_{\mathcal A}\cap \overline B_{1-\epsilon} \mid \exists \rho>0 \text{ with } \spt T_{\mathcal A}\cap B_\rho(x)=\spt T_{\mathcal E}\cap B_\rho(x) \bigg \}.\]
Clearly $\mathcal Z_\epsilon \neq \emptyset$ since $0 \in \mathcal Z_\epsilon$ and it is an open subset of $\spt T_{\mathcal A}\cap \overline B_{1-\epsilon}.$ Note that by transversality at $\partial B_1$ we have that $\spt T_{\mathcal A}\cap  \overline B_{1-\epsilon}$ is connected for small enough $\epsilon.$ We show that $\mathcal Z_\epsilon$ is also closed.

Let $x_i \in \mathcal Z_\epsilon $ with $x_i \to x_0.$ We assume that $x_0 \not= 0.$ Clearly $x_0 \in \spt T_{\mathcal A} \cap \spt T_{\mathcal E} \cap \overline B_{1-\epsilon}.$ Thus it is a regular point of $\spt T_{\mathcal A}$. The variational equations are valid around $x_0$ for both. Blowing up the associated varifolds around $x_0$, we get a stationary cone lying in a half--space. From \cite{Simon's notes}, we conclude that the limit cone of $V_{\mathcal E}$ at $x_0$ is a plane (possibly with multiplicity). Since we have no touching singularities for $V_{\mathcal E}$, the sheeting theorem of \cite{BW1} implies that $x_0 $ is a regular point of $\spt T_{\mathcal E}$. Thus the standard maximum principle for CMC hypersurfaces implies that $\spt T_{\mathcal A}$ and $\spt T_{\mathcal E}$ coincide around $x_0$, hence $x_0 \in \mathcal Z_\epsilon.$ This shows that $\mathcal Z_\epsilon$ is closed. Letting $\epsilon \to 0$, we conclude that $\spt T_{\mathcal A}\cap  B_1\subset\spt T_{\mathcal E}\cap   B_1.$  

 Taking closures and using Lemma \ref{var part I} and Lemma \ref{avoidance}, we get that $\spt \partial T_{\mathcal A}\subset \spt \partial T_{\mathcal E}$, which implies that the boundaries must agree since they are connected, embedded submanifolds of the same dimension. Thus the set $F = (A\cap B_1)\cup (E\setminus B_1)$ is a competitor for the functional $J_\lambda$ and satisfies $\Per_{\partial B_1}(F)=\Per_{\partial B_1}(E)$. Since $\Per_{B_1}(A)\leq \Per_{B_1}(E)$, we get that
\[0\leq J_\lambda (F)-J_\lambda(E) \leq -\lambda |(A\setminus E) \cap B_1|, \] 
which implies that $A\cap B_1=E\cap B_1.$ 
\end{proof}

\section{Smooth approximation} 
\noindent
In this section we prove the smooth approximation theorem. As before, let $A \subset B_2 \subset  \mathbb R^8$ be a set of finite perimeter that minimizes $J_\lambda$ on one side, for $0<\lambda<7$, and is the unique minimizer of $J_\lambda$ with respect to the boundary condition $A \setminus (A\cap B_1)$. Furthermore, $\partial A$ is smoothly embedded except at the origin,  the associated varifold $V_{\mathcal A}$ satisfies \eqref{variational eq} and $\partial A \cap \partial B_1$ is connected and a transverse intersection. Finally, we assume that a varifold tangent of $V_{\mathcal A}$ at the origin is a  non--flat regular cone $C$ of multiplicity--1 that locally minimizes the perimeter functional on one side. Thus $\mathbb R ^8 \setminus C$ consists of two connected components $D_\pm$. For the associated current, we write $C= \partial \llbracket D\rrbracket$, where $D=D_+$. We show the following: 

\begin{thm}\label{smooth app1}(smooth approximation in $\mathbb R ^8$) There exists a sequence of sets of finite perimeter $A_j$ that minimize $J_\lambda$ on one side and $A_j \to A.$ Furthermore, there exists a ball $B_r \subset B_1$ such that $\partial A_j\cap B_r$ is a smoothly embedded CMC hypersurface that converges smoothly to $\partial A \cap B_r$ on any $K \subset \subset B_r\setminus \{0\}.$  
\end{thm}

\begin{proof} We divide the proof in steps. 

\noindent
\emph{Step 1.} We use the same perturbation as in the proof of \cite[Theorem 5]{BL} to perturb the boundary condition inwards. Let $d(\cdot):=d(\cdot, \partial A)$ be the signed distance function of $\partial A$, which is positive in $A$, and consider $T_\delta$ a tubular neighborhood of $\partial A \cap (B_{\frac{3}{2}}\setminus B_{\frac{1}{2}})$ where the signed distance is smooth. Let $\eta\in C^1_c(B_2)$ such that $0\leq \eta \leq 1$, $\eta \equiv 1$ in $T_{\delta/2} \cap (B_{\frac{5}{4}}\setminus B_{\frac{3}{4}})$ and $\spt \eta \subset T_\delta\cap (B_{\frac{3}{2}}\setminus B_{\frac{1}{2}})$ and define the vector field $X = \eta \nabla d.$ Then $X$ is a smooth extension of the unit normal $\nu_{ A}$ near $\partial A \cap \partial B_1$. For $t \in [0, \epsilon)$, we consider the flow $\phi_t$ of $X.$ Then $\phi_t(A) \subset A$, $\phi_t (\partial A) \cap \partial B_1$ is disjoint from $\partial A \cap \partial B_1$ and $\phi_t (A)\cap (B_{\frac{5}{4}}\setminus B_{\frac{3}{4}})$ is strictly contained in $A \cap (B_{\frac{5}{4}}\setminus B_{\frac{3}{4}})$ for $0<t<\epsilon.$

Choose a sequence $t_j \searrow 0$ and consider the minimization problem 
\begin{equation}
\label{perturbed min problem}
\inf\bigg\{J_\lambda(E)\mid E\subset B_2 ,  E \setminus (A\cap B_1) = \phi_{t_j}(A)\setminus (A\cap B_1) \bigg \}.
\end{equation}
From Proposition \ref{existence}, there are sets of finite perimeter $A_{j}$ that minimize \eqref{perturbed min problem}. 

From the results of Section \ref{2.3}, the associated varifolds $V_{\mathcal A_j}$ satisfy the variational equations  \eqref{variational eq} away from their current boundary $T_{b,j}$. Furthermore, from Theorem \ref{reg theorem 2} we have that $\sing_T V_{\mathcal A_j}=\emptyset$ and $\sing V_{\mathcal A_j}$ consists of isolated points.  

\medskip
\noindent 
\emph{Step 2.} We show for the associated varifolds $V_{\mathcal A_j}, V_{\mathcal A}$ that 
\[V_{\mathcal A_j} \xrightarrow[j\to \infty]{} V_{\mathcal A},\]
in the varifold topology. 

By construction, we have that $\phi_{t_j}(A) \xrightarrow[j\to \infty]{} A$. Using the area formula we have that $J_\lambda(\phi_{t_j}(A)) \to J_\lambda(A)$ as $j \to \infty.$ From the minimization property we have that $J_\lambda(A_{j}) \leq J_\lambda(\phi_{t_j}(A))$ and so $\sup_j\Per_{B_2}(A_{j}) < \infty$. From the standard compactness for sets of finite perimeter, see \cite{Maggi}, we conclude that $A_{j} \to \tilde A$ with $\tilde A \setminus (A\cap B_1) = A \setminus (A\cap B_1)$ and $\tilde A \subset A.$ From the lower semi--continuity we obtain that $J_\lambda(\tilde A) \leq J_\lambda (A)$. Thus $\tilde A$ is a one--sided minimizer with respect to the boundary condition $A \setminus (A\cap B_1)$ and so by the uniqueness assumption $\tilde A = A.$

In particular, we have shown that $A_{j} \to A$ as sets of finite perimeter. Finally, the minimization property combined with the area formula gives that 
\[\Per_{B_2}(A_{j})-\lambda |A_{j}| \leq \Per_{B_2}(A)(1 + O(t_j)) -\lambda |\phi_{t_j}(A)|.\] 
Combining $A_j \to A$ and the preceding inequality with the lower semi--continuity of the perimeter functional we conclude, after passing to a subsequence, that 
\[\Per_{B_2}(A_{j}) \to \Per_{B_2}(A).\]
This is enough to get the varifold convergence in $B_1.$  

\medskip
\noindent 
\emph{Step 3.} It remains to show that there exists $r\in (0,1)$  such that $\partial A_{j} \cap B_r$ is smooth for all large enough $j$'s. Assume to the contrary that there exists a sequence $r_j \searrow 0$ such that $\partial A_{j} \cap B_{r_j}$ is not smooth. Thus there exists a sequence of points $p_j \in \sing\partial A_{j} \cap B_{r_j}$ and $p_j \to 0.$  Then $p_j \neq 0$ for all $j$ since otherwise $0 \in \partial A \cap \partial A_{j}$ and the maximum principle, Corollary \ref{singular max pri}, implies that $\partial A_j \cap \overline B_1 = \partial A \cap \overline B_1$ which is not possible. 

Define the sequence $O_j :=\eta_{0, |p_j|}(A_{j})$ and let $\ V_{\mathcal O _j} $ be the associated varifold. For any $R>0$, the variational equations \eqref{variational eq}, Lemma \ref{stability} and Theorem \ref{reg theorem 2} imply that the varifold $ V_{\mathcal O _j}$ satisfies all the assumptions of the compactness theorem of \cite{BW1} in $B_R.$ Thus $ V_{\mathcal O_j} \to  V_\infty,$ where $V_\infty$ is an integral varifold that is stationary, stable and has no classical singularities. By \cite{Wick2}, $\sing V_\infty$ consists of isolated points. Finally, note that $\spt V_\infty \subset \overline D$. (Here we use the uniqueness of the varifold tangent $C$).

We show that the regular part of $V_\infty$ has multiplicity--1 and is connected. From the varifold convergence of $V_{\mathcal A_j}$ to $ V_{\mathcal A}$ and the monotonicity of the mass ratio we conclude that 
\begin{equation}
\label{blow-down}
\Theta^7_\infty(V_\infty)=\lim_{R\to \infty}\dfrac{\mu_{V_\infty}(B_R)}{\omega_7 R^7 }\leq \Theta^7(C,0). 
\end{equation}
Let $\Sigma_\alpha$ be the connected components of $\reg V_\infty$. Thus we may decompose 
\[V_\infty = \sum_\alpha k_\alpha |\Sigma_\alpha|, \]
where $k_\alpha \in \mathbb N$ and $|\Sigma_\alpha|$ is the associated varifold to $\Sigma_\alpha.$ Note that since the singular set of $V_\infty$ consists of isolated points  each $|\Sigma_\alpha|$ is stationary (see \cite{Wick max} for a similar argument). 

From \eqref{blow-down} we have that for each varifold $|\Sigma_\alpha|$ the blow--down limit exists. In particular, from Allard's compactness theorem, we get a varifold tangent cone $C_\alpha$ at infinity that is stationary. Note that the inclusion $\spt V_\infty \subset \overline D$ implies that $\spt C_\alpha \subset \overline D$. The maximum principle of \cite{Wick max} implies that $\spt C_\alpha = \spt C$ and from the constancy theorem we get that $C_\alpha = l C$, where $l \in \mathbb N.$ 

Furthermore, from \eqref{blow-down} we have that 
\[l \Theta^7(C,0)= \Theta^7(C_\alpha, 0)= \Theta_\infty (|\Sigma_\alpha|)\leq k_\alpha^{-1} \Theta^7_\infty(V_\infty)\leq k_\alpha^{-1}\Theta^7(C,0).\]
In particular, $l=k_\alpha=1.$ This shows that each connected component of the regular part has multiplicity--1. 

Similarly, if we have two connected components we would get 
\[2\Theta^7(C,0)\leq \Theta^7_\infty(|\Sigma_1|)+\Theta^7_\infty(|\Sigma_2|)\leq \Theta^7_\infty(V_\infty)\leq \Theta^7(C,0),\]
which is not possible. This shows that $\reg V_\infty$ is also connected. Since the singular set of $V_\infty$ consists of isolated points we have that $\spt V_\infty = \partial O$ where $O$ is a set of locally finite perimeter, with $O \subset D$ and from the constancy theorem we conclude that $V_\infty = \underbar v(\partial O , 1). $ 

Note that the multiplicity--1 convergence implies, from \cite{White2}, that for the associated currents we have that $T_{\mathcal O_j}\rightharpoonup \partial \llbracket O \rrbracket$ and that $O_j \to O$ as sets of finite perimeter. 

\medskip
\noindent
\emph{Step 4.} We now show that the limit set $O$ locally minimizes the perimeter functional on one side. Given $x\in \partial O$ we want to find $r>0$ such that $\Per_{B_{2r}(x)}(O)\leq \Per_{B_{2r}(x)}(F)$ for any $F \subset B_{2r}(x)$ with $F\setminus (B_r(x) \cap O)= (O\cap B_{2r}(x))\setminus (B_r(x)\cap O).$ We will abbreviate the ball by $B_r.$

Define the sets $F_j = (F\cap O_j \cap B_r)\cup (B_{2r}\setminus B_r) \cap O_j. $ Note that $O_j\cap B_{2r}$ minimizes $J_{\lambda_j}$ on one side, where $\lambda_j = \lambda |p_j|$.  Since $F_j \subset O_j$ and satisfies the boundary condition, the minimization condition gives that 
\[\Per_{B_r}(O_j) +\Per_{\partial B_r}(O_j) - \lambda_j|O_j \cap B_r| \leq \Per_{B_r}(F_j) +\] \[+\Per_{\partial B_r}(F_j)-\lambda_j |F_j \cap B_r|.\]
Using the formula in \cite[Theorem 16.3]{Maggi} we have that $\Per_{B_r}(F_j) = \Per_{B_r}(F\cap O_j).$ To evaluate the term $\Per_{\partial B_r}(F_j)$ we choose $r>0$ such that $\partial O \cap \partial B_r$ is a transverse intersection (such a choice is possible by Sard's theorem since $\partial O$ is embedded except at isolated singular points). Then from the established multiplicity--1 varifold convergence of $V_{\mathcal O_j}$ to $V_{\mathcal O}$ we can apply Allard's interior regularity theorem \cite{Allard1}, to get that for all large $j$ the intersection $\partial O_j \cap \partial B_r$ is also transverse. 

Using the formula in \cite[Theorem 16.3]{Maggi} and the transversality condition we get that 
\[\Per_{\partial B_r}(F_j) = \mathcal H^7(O_j \cap \partial B_r\cap F^{(0)})+ \mathcal H^7(O_j \cap \partial B_r\cap \partial^* F \cap \{\nu_F= \nu_{B_r^c}\}),\] 
where $F^{(0)}$ is the set of points of $F$ of zero $8$--dimensional density. 
Using \cite[Lemma 12.22]{Maggi} for the term $\Per_{B_r}(F\cap O_j)$ we get that 
\[\Per_{B_r}(O_j \cup F) -\lambda_j |O_j \cap B_r| \leq \Per_{B_r}(F) +\mathcal H^7(\partial B_r\cap O_j \cap F^{(0)})+\] \[+\mathcal H^7(O_j \cap \partial B_r\cap \partial^* F \cap \{\nu_F= \nu_{B_r^c}\})-\lambda_j |F\cap O_j\cap B_r|.\]

The established set and multiplicity--1 varifold convergence of the previous step allow us to use \cite{Lahti} to get that
\[\chi_{O_j}\to \chi_O \text{ for } \mathcal H^7 \text{ --a.e. point in } \mathbb R ^8 \setminus \partial O.\]
By the transversality assumption we have that convergence holds for $\mathcal H^7$--a.e. point in $\partial B_r.$ Thus from the dominated convergence theorem we conclude that for $j \to \infty$
\[\mathcal H ^7(\partial B_r \cap O_j \cap F^{(0)}) \to \mathcal H^7(\partial B_r \cap O \cap F^{(0)})\]
\[\mathcal H^7(O_j \cap \partial B_r\cap \partial^* F \cap \{\nu_F= \nu_{B_r^c}\})\to \mathcal H^7(O\cap \partial B_r\cap \partial^* F \cap \{\nu_F= \nu_{B_r^c}\} ).\]

Combining $O_j \cup F \to O$, the vanishing of the volume terms and the lower semi--continuity of the perimeter functional, we obtain 

\[\Per_{B_r}(O) \leq \Per_{B_r}(F) + \mathcal H^7( O\cap \partial B_r \cap F^{(0)})+  \]
\[+\mathcal H^7(O\cap \partial B_r\cap \partial^* F \cap \{\nu_F= \nu_{B_r^c}\} ).\]
Using again the formula of \cite[Theorem 16.3]{Maggi} and the transversality condition we get that the sum of the last two terms equals to $\Per_{\partial B_r}(F\cap O \cap B_r \cup (B_{2r}\setminus B_r)\cap O) = \Per_{\partial B_r}(F). $ 

In total, using again the transversality and after adding and subtracting the boundary term $\Per_{B_{2r}\setminus \overline B_r}(O)$ we get that 
\[\Per_{B_{2r}}(O)\leq \Per_{B_{2r}}(F),\]
thus showing the one--sided minimization property. 

\medskip
\noindent
\emph{Step 5.} Since $O \subset D$ and $O, D$ locally minimize the perimeter functional on one side, we may apply the Hardt--Simon foliation \cite{HS}, extended by \cite{Liu}, to get that $\partial O$ is either the cone itself or a leaf of the foliation. In any case $\partial O $ is smooth in a neighborhood of $S^7$. 

We can now obtain a contradiction using Allard's regularity theorem as in \cite{BL} or \cite{HS}. We briefly recall the details. Since up to a subsequence, $p_j|p_j|^{-1}\to \xi \in S^7,$ upper semi--continuity of the density implies that 
\[\limsup_j\Theta^7(V_{\mathcal O_j}, p_j|p_j|^{-1})\leq \Theta^7(\partial O , \xi) =1.\]
Using the scale invariance of the density we obtain that for large $j$ 
\[\Theta^7(V_{\mathcal A_j}, p_j)< 1 +\epsilon,\]
where $\epsilon$ is the universal constant from Allard's interior regularity theorem \cite{Allard1}. But then it follows that $p_j \in \reg V_{\mathcal A_j}$, contrary to our assumption. 

The stated smooth convergence follows from Step 2, Allard's interior regularity theorem and standard elliptic regularity. 
\end{proof}

We explain why Theorem \ref{smooth app1} is enough to get Theorem \ref{main theorem}. Let $A \subset \mathbb R^8$ be a set of locally finite perimeter that minimizes $J_\lambda$ with respect to the boundary condition $(B_{2R}\setminus B_R)\cap A$ for some $R>0$ and $\lambda>0.$ Furthermore the associated varifold $V_{\mathcal A}$ satisfies the variational equations \eqref{variational eq}. Then due to Theorem \ref{reg theorem 1}  we may repeat the same analytic continuation argument of \cite[Lemma 5.1]{BL}.  This shows that if $R'< R$ then $A$ is the unique one sided minimizer with respect to the boundary condition $(B_{2R'}\setminus B_{R'})\cap A.$ Thus we may choose a small enough $R'>0$ so that $\lambda < \frac{7}{R'}$ and connectedness of the boundary also holds.  Scaling the set $A$ by $R'$ we get the working conditions of Theorem \ref{smooth app1}.

\begin{rem}\label{satisfying the assumptions} We mention examples of one--sided minimizers that satisfy the assumptions of Theorem \ref{main theorem}. Consider the cone  $C:= C_{1,5}=\{(x,y)\in \mathbb R ^2\times \mathbb R^6 \mid 5|x|^2=|y|^2 \}$ which minimizes the perimeter functional on one side, see \cite{Dav}. Then from \cite{CHS}, extended to the CMC setting in \cite{Irv}, there exist CMC hypersurfaces asymptotic to the cone $C.$ Using arguments similar to those in \cite{HS} as in \cite{Lesk} (see also the recent work of \cite{Kobe} on singular isoperimetric boundaries), one shows that these hypersurfaces minimize the prescribed mean curvature functional on one side and in particular satisfy all the assumptions of Theorem \ref{main theorem}. 
    
\end{rem} 

\begin{rem}\label{on the sign of lambda} We only considered the case where the prescribed scalar mean curvature is positive. If $\lambda<0$, all the arguments of Section \ref{2.3} extend in a straightforward manner under the assumption that $|\lambda|< n.$ Similarly the results of Section \ref{3} follow under the assumption that $|\lambda|<7.$  The only argument that needs to be adjusted is the final inequality in Corollary \ref{singular max pri}. In particular, once we have established 
\[0\leq J_\lambda (F)-J_\lambda(E) \leq -\Per_{B_1}(F\setminus E)-\lambda |F\setminus E|\]
then, since $\Per_{\mathbb R^8}(F\setminus E) = \Per_{B_1}(F\setminus E)$, a standard application of the divergence theorem, see \cite[Remark 12.36]{Maggi} gives that the latter term is 
\[\leq -(8-|\lambda|)|(F\setminus E)\cap B_1|\] 
and the assumption that $|\lambda| < 7$ implies the final result. In conclusion Theorem \ref{smooth app1} holds for $|\lambda |<7$.
    
\end{rem}

\begin{rem}\label{any dimension} Theorem \ref{smooth app1} extends to the case that $A \subset B_2 \subset  \mathbb R^{n+1}$ for $n\geq 7.$ We assume the same hypotheses as in Theorem \ref{smooth app1} with $\mathbb R^{n+1}$ instead of $\mathbb R ^8$ and $\lambda <n.$ The proof of Theorem \ref{mult--1} extends using the same arguments. One needs to substitute the transversality condition by choosing radii so that the spherical slices of the varifolds have zero $n$--dimensional Hausdorff measure.  Furthermore, from regularity theory of \cite{BW1} and \cite{Wick2} the singular maximum principle argument \cite{Wick max} still applies. This establishes Theorem \ref{mult--1} and hence Corollary \ref{singular max pri} follows directly. In the proof of Theorem \ref{smooth app1} one uses the regularity theory of \cite{BW1} and \cite{Wick2}. Then the connected components $|\Sigma_\alpha|$ will be stationary since the argument of \cite{Wick max} holds even if the singularities are not isolated. As before, one needs to substitute transversality with spherical slices of zero $n$--dimensional Hausdorff measure and use the exact same argument to get the smooth approximation. Thus Theorem \ref{smooth app1} and hence Theorem \ref{main theorem} hold for general $n \geq 7$.
\end{rem}

\end{document}